\input amstex
\input amsppt.sty
\magnification=\magstep1
\hsize=32truecc
\vsize=22.2truecm
\baselineskip=16truept
\NoBlackBoxes
\TagsOnRight \pageno=1 \nologo
\def\Z{\Bbb Z}
\def\N{\Bbb N}

\def\Q{\Bbb Q}

\def\l{\left}
\def\r{\right}
\def\bg{\bigg}
\def\({\bg(}
\def\[{\bg\lfloor}
\def\){\bg)}
\def\]{\bg\rfloor}
\def\t{\text}
\def\f{\frac}

\def\bi{\binom}
\def\eq{\equiv}

\def\ls{\leqslant}
\def\gs{\geqslant}
\def\mo{\roman{mod}}

\def\Remark{\medskip\noindent{\it  Remark}}

\hbox {{\tt arXiv:1102.5649}}
\bigskip
\topmatter
\title List of conjectural series for powers of $\pi$ and other constants\endtitle
\author Zhi-Wei Sun\endauthor
\leftheadtext{Zhi-Wei Sun} \rightheadtext{Series for
powers of $\pi$ and other constants}
\affil Department of Mathematics, Nanjing University\\
 Nanjing 210093, People's Republic of China
  \\  zwsun\@nju.edu.cn
  \\ {\tt http://math.nju.edu.cn/$\sim$zwsun}
\endaffil
\abstract Here I give the full list of my conjectures on series for
powers of $\pi$ and other important constants scattered in some of
my public papers or my private diaries. The list contains 234 reasonable conjectural series. On the list there are 178 reasonable series for
$\pi^{-1}$, four series for $\pi^2$, two series for $\pi^{-2}$, four series for $\pi^4$, two series for $\pi^5$, three series for $\pi^6$, seven
series for $\zeta(3)$, one series for $\pi\zeta(3)$, two series for $\pi^2\zeta(3)$, one series for $\zeta(3)^2$, three series involving both $\zeta(3)^2$ and $\pi^6$, one series for $\zeta(5)$, three series involving both $\zeta(5)$ and $\zeta(2)\zeta(3)$, two series involving both $\pi\zeta(5)$ and $\pi^3\zeta(3)$, three series involving $\zeta(7)$, three series for $K=L(2,(\frac{\cdot}{3}))$, one series for the Catalan constant $G$, two series for $\pi G$, one series involving both $\pi^3G$ and $\pi^2\zeta(3)$, two series for $\pi K$, two series involving $L=L(4,(\frac{\cdot}3))$, three series involving $\beta(4)=L(4,(\frac{-4}{\cdot}))$, and four series for $\pi^2\log a$ with $a=2,3,(\sqrt5+1)/2$. The
code of a conjectural series is underlined if and only if a complete proof of the identity is available.
\endabstract
\thanks \copyright\ Copyright is owned by the author Zhi-Wei Sun.
The list on the author's homepage has been linked to {\tt Number Theory Web} since
Feb. 18, 2011.
\endthanks
\endtopmatter
\document

\heading{1. Series for various constants other than $1/\pi$}\endheading

\proclaim{Conjecture 1} {\rm (i) (Z. W. Sun [S11, Conj. 1.4])} We have
$$\align\sum_{k=1}^\infty\f{(10k-3)8^k}{k^3\bi{2k}k^2\bi{3k}k}=&\f{\pi^2}2,\tag $\underline{1.1}$
\\\sum_{k=1}^\infty\f{(11k-3)64^k}{k^3\bi{2k}k^2\bi{3k}k}=&8\pi^2,\tag $\underline {1.2}$
\\\sum_{k=1}^\infty\f{(35k-8)81^k}{k^3\bi{2k}k^2\bi{4k}{2k}}=&12\pi^2,\tag $\underline{1.3}$
\endalign$$
$$\align\sum_{k=1}^\infty\f{(15k-4)(-27)^{k-1}}{k^3\bi{2k}k^2\bi{3k}k}=&K,\tag$\underline {1.4}$
\\\sum_{k=1}^\infty\f{(5k-1)(-144)^k}{k^3\bi{2k}k^2\bi{4k}{2k}}=&-\f{45}2K,\tag $\underline{1.5}$
\endalign$$ where
$$\align K:=&L\l(2,\l(\f {\cdot}3\r)\r)=\sum_{k=1}^\infty\f{(\f k3)}{k^2}
\\=&0.781302412896486296867187429624\ldots
\endalign$$ with $(-)$ the Legendre symbol.
\medskip

{\rm (ii) (Z. W. Sun [S13a, Conj. 8])} We have
$$\sum_{k=1}^\infty\f{(28k^2-18k+3)(-64)^k}{k^5\bi{2k}k^4\bi{3k}k}=-14\zeta(3),\tag1.6$$
where $\zeta(3):=\sum_{n=1}^\infty1/n^3.$

{\rm (iii) (Z. W. Sun [S11, Conj. 1.4])} We have
$$\sum_{n=0}^\infty\f{18n^2+7n+1}{(-128)^n}\bi{2n}n^2\sum_{k=0}^n\bi{-1/4}k^2\bi{-3/4}{n-k}^2=\f{4\sqrt2}{\pi^2}\tag1.7$$
and
$$\sum_{n=0}^\infty\f{40n^2+26n+5}{(-256)^n}\bi{2n}n^2\sum_{k=0}^n\bi{n}k^2\bi{2k}k\bi{2(n-k)}{n-k}
=\f{24}{\pi^2}.\tag1.8$$

{\rm (iv) (Z. W. Sun [S14e, Conj. 1.1])} We have
$$\sum_{k=1}^\infty\f{48^k}{k(2k-1)\bi{4k}{2k}\bi{2k}k}=\f{15}2K,\tag1.9$$
where $K$ is as in part {\rm (i)}.

{\rm (v) (Z. W. Sun [S14e, Conj. 3.1 and 3.2])}
For $n=1,2,3,\ldots$ let $H_n$ denote the harmonic number $\sum_{k=1}^n1/k$. Then
$$\align\sum_{k=1}^\infty\f{H_{2k}+2/(3k)}{k^2\bi{2k}k}=&\zeta(3),\tag1.10
\\\sum_{k=1}^\infty\f{H_{2k}+2H_k}{k^2\bi{2k}k}=&\f53\zeta(3),\tag1.11
\\\sum_{k=1}^\infty\f{H_{2k}+17H_k}{k^2\bi{2k}k}=&\f52\sqrt3\ \pi K.\tag1.12
\endalign$$
Also,
$$\align\sum_{k=1}^\infty\f{2^k}{k^2\bi{2k}k}\l(H_{\lfloor k/2\rfloor}-(-1)^k\f 2k\r)=&\f 74\zeta(3),\tag1.13
\\\sum_{k=1}^\infty\f{2^k}{k^2\bi{2k}k}\l(2H_{2k}-3H_k+\f2k\r)=&\f 74\zeta(3),\tag1.14
\\\sum_{k=1}^\infty\f{2^k}{k^2\bi{2k}k}\l(6H_{2k}-11H_k+\f 8k\r)=&2\pi G,\tag1.15
\\\sum_{k=1}^\infty\f{2^k}{k^2\bi{2k}k}\l(2H_{2k}-7H_k+\f2k\r)=&-\f{\pi^2}2\log2,\tag1.16
\\\sum_{k=1}^\infty\f{3^k}{k^2\bi{2k}k}\l(6H_{2k}-8H_k+\f5k\r)=&\f{26}3\zeta(3),\tag1.17
\\\sum_{k=1}^\infty\f{3^k}{k^2\bi{2k}k}\l(6H_{2k}-10H_k+\f7k\r)=&2\sqrt3\,\pi K,\tag1.18
\\\sum_{k=1}^\infty\f{3^k}{k^2\bi{2k}k}\l(H_k+\f1{2k}\r)=&\f{\pi^2}3\log3,\tag1.19
\\\sum_{k=1}^\infty\f{L_{2k}}{k^2\bi{2k}k}\l(\f1k+\f1{k+1}+\cdots+\f1{2k}\r)=&\f{41\zeta(3)+4\pi^2\log\phi}{25},\tag1.20
\\\sum_{k=1}^\infty\f{v_k}{k^2\bi{2k}k}\l(\f1k+\f1{k+1}+\cdots+\f1{2k}\r)=&\f{124\zeta(3)+\pi^2\log(5^5\phi^6)}{50},\tag1.21
\endalign$$
where $G$ denotes the Catalan constant $\sum_{k=0}^\infty(-1)^k/(2k+1)^2$, $\phi$ stands for the famous golden ratio $(\sqrt5+1)/2$, the Lucas numbers $L_0,L_1,L_2,\ldots$
are given by 
$$L_0=2,\ L_1=1,\ \t{and}\ L_{n+1}=L_n+L_{n-1}\ \t{for}\ n=1,2,3,\ldots,$$  and $v_0,v_1,v_2,\ldots$ 
are defined by
$$v_0=2,\ v_1=5,\ \t{and}\ v_{n+1}=5(v_n-v_{n-1})\ \t{for}\ n=1,2,3,\ldots.$$

{\rm (vi) (Z. W. Sun [S14e, Conj. 3.3 and 3.4])}
$$\align\sum_{k=1}^\infty(-1)^{k-1}\f{10H_k-3/k}{k^3\bi{2k}k}=&\f{\pi^4}{30},\tag1.22
\\\sum_{k=1}^\infty(-1)^{k-1}\f{H_{2k}+4H_k}{k^3\bi{2k}k}=&\f2{75}\pi^4,\tag1.23
\\\sum_{k=1}^\infty\f{H_{2k}-H_k+2/k}{k^4\bi{2k}k}=&\f{11}9\zeta(5),\tag1.24
\endalign$$
$$\align\sum_{k=1}^\infty\f{3H_{2k}-102H_k+28/k}{k^4\bi{2k}k}=&-\f{55}{18}\pi^2\zeta(3),\tag1.25
\\\sum_{k=1}^\infty\f{97H_{2k}-163H_k+227/k}{k^4\bi{2k}k}=&\f{165}8\sqrt3\pi L,\tag1.26
\endalign$$
where
$$L:=L\l(4,\l(\f{\cdot}3\r)\r)=\sum_{k=1}^\infty\f{(\f k3)}{k^4}.$$
We also have
$$\align\sum_{k=0}^\infty\f{\bi{2k}k}{(2k+1)16^k}\l(3H_{2k+1}+\f4{2k+1}\r)=&8G,\tag1.27
\\\sum_{k=0}^\infty\f{\bi{2k}k}{(2k+1)^2(-16)^k}\l(5H_{2k+1}+\f{12}{2k+1}\r)=&14\zeta(3),\tag1.28
\\\sum_{k=0}^\infty\f{\bi{2k}k}{(2k+1)^316^k}\l(9H_{2k+1}+\f{32}{2k+1}\r)=&40\beta(4)+\f5{12}\pi\zeta(3),\tag1.29
\endalign$$
where $$\beta(4)=L\l(4,\l(\f{-4}{\cdot}\r)\r)=\sum_{k=0}^\infty\f{(-1)^k}{(2k+1)^4}.$$

{\rm (vii) (Z. W. Sun [S14e, Conj. 4.1 and 4.2]) Let $H_k^{(m)}$ denote $\sum_{0<j\ls k}1/j^m$. Then
$$\sum_{k=1}^\infty\f{6H_{\lfloor k/2\rfloor}^{(2)}-(-1)^k/k^2}{k^2\bi{2k}k}=\f{13}{1620}\pi^4.\tag1.30$$
Also,
$$\align\sum_{k=1}^\infty\f{H_k^{(3)}}{k^2\bi{2k}k}=&\f{\zeta(5)+2\zeta(2)\zeta(3)}9,\tag1.31
\\\sum_{k=1}^\infty\f{(-1)^k}{k^3\bi{2k}k}\(10\sum_{j=1}^k\f{(-1)^j}{j^2}-\f{(-1)^k}{k^2}\)=&\f{29\zeta(5)-2\zeta(2)\zeta(3)}6,\tag1.32
\\\sum_{k=1}^\infty\f1{k^2\bi{2k}k}\(24\sum_{j=1}^k\f{(-1)^j}{j^3}-17\f{(-1)^k}{k^3}\)=&7\zeta(5)-6\zeta(2)\zeta(3),\tag1.33
\endalign$$
$$\align\sum_{k=1}^\infty(-1)^{k-1}\f{H_k^{(3)}+1/(5k^3)}{k^3\bi{2k}k}=&\f25\zeta(3)^2,\tag1.34
\\\sum_{k=1}^\infty\f{H_{k-1}^{(2)}-1/k^2}{k^4\bi{2k}k}=&-\f{313\pi^6}{612360},\tag1.35
\\\sum_{k=1}^\infty\f{3H_k^{(4)}-1/k^4}{k^2\bi{2k}k}=&\f{163\pi^6}{136080},\tag1.36
\endalign$$
$$\align\sum_{k=1}^\infty\f1{k^4\bi{2k}k}\(72\sum_{j=1}^k\f{(-1)^j}{j^2}-\f{(-1)^k}{k^2}\)=&-\f{31}{1134}\pi^6-\f{34}5\zeta(3)^2,\tag1.37
\\\sum_{k=1}^\infty\f1{k^2\bi{2k}k}\(8\sum_{j=1}^k\f{(-1)^j}{j^4}+\f{(-1)^k}{k^4}\)=&-\f{97}{34020}\pi^6-\f{22}{15}\zeta(3)^2,\tag1.38
\\\sum_{k=1}^\infty\f{(-1)^k}{k^3\bi{2k}k}\(40\sum_{0<j<k}\f{(-1)^j}{j^3}-7\f{(-1)^k}{k^3}\)=&-\f{367}{27216}\pi^6+6\zeta(3)^2.\tag1.39
\endalign$$

{\rm (viii)} (Z. W. Sun [S14e, Conj. 4.3]) We have
$$\align\sum_{k=1}^\infty\f{33H_k^{(5)}+4/k^5}{k^2\bi{2k}k}=&-\f{45}8\zeta(7)+\f{13}3\zeta(2)\zeta(5)+\f{85}6\zeta(3)\zeta(4),\tag1.40
\\\sum_{k=1}^\infty\f{33H_k^{(3)}+8/k^3}{k^4\bi{2k}k}
=&-\f{259}{24}\zeta(7)-\f{98}9\zeta(2)\zeta(5)+\f{697}{18}\zeta(3)\zeta(4),\tag1.41
\endalign$$
and
$$\sum_{k=1}^\infty\f{(-1)^k}{k^3\bi{2k}k}\(110\sum_{j=1}^k\f{(-1)^j}{j^4}+29\f{(-1)^k}{k^4}\)
=\f{223}{24}\zeta(7)-\f{301}6\zeta(2)\zeta(5)+\f{221}2\zeta(3)\zeta(4).\tag1.42$$

{\rm (ix)} (Z. W. Sun [S14e, Conj. 5.1-5.3]} We have
$$\align\sum_{k=0}^\infty\f{\bi{2k}k}{(2k+1)16^k}\sum_{j=0}^k\f1{(2j+1)^3}=&\f 5{18}\pi\zeta(3),\tag1.43
\\\sum_{k=0}^\infty\f{\bi{2k}k}{(2k+1)^2(-16)^k}\sum_{j=0}^k\f{(-1)^j}{(2j+1)^2}=&\f{\pi^2G}{10}+\f{\pi\zeta(3)}{240}+\f{27\sqrt3}{640}L,\tag1.44
\endalign$$
$$\align\sum_{k=0}^\infty\f{\bi{2k}k}{(2k+1)8^k}\(\sum_{j=0}^k\f{(-1)^j}{2j+1}-2\f{(-1)^k}{2k+1}\)=&-\f{\sqrt2}{16}\pi^2,\tag1.45
\\\sum_{k=0}^\infty\f{\bi{2k}k}{(2k+1)16^k}\(12\sum_{j=0}^k\f{(-1)^j}{(2j+1)^2}-\f{(-1)^k}{(2k+1)^2}\)=&4\pi G,\tag1.46
\\\sum_{k=0}^\infty\f{\bi{2k}k}{(2k+1)^2(-16)^k}\(5\sum_{j=0}^k\f1{(2j+1)^3}+\f1{(2k+1)^3}\)=&\f{\pi^2}2\zeta(3),\tag1.47
\\\sum_{k=0}^\infty\f{\bi{2k}k}{(2k+1)16^k}\(24\sum_{j=0}^k\f{(-1)^j}{(2j+1)^3}-17\f{(-1)^k}{(2k+1)^3}\)=&\f{\pi^4}{12},\tag1.48
\\\sum_{k=0}^\infty\f{\bi{2k}k}{(2k+1)^2(-16)^k}\(40\sum_{j=0}^k\f{(-1)^j}{(2j+1)^3}-47\f{(-1)^k}{(2k+1)^3}\)=&-\f{85\pi^5}{3456},\tag1.49
\\\sum_{k=0}^\infty\f{\bi{2k}k}{(2k+1)16^k}\(3\sum_{j=0}^k\f1{(2j+1)^4}-\f1{(2k+1)^4}\)=&\f{121\pi^5}{17280},\tag1.50
\\\sum_{k=0}^\infty\f{\bi{2k}k}{(2k+1)^2(-16)^k}\(5\sum_{j=0}^k\f1{(2j+1)^4}-\f 4{(2k+1)^4}\)=&\f{7\pi^6}{7200}.\tag1.51
\endalign$$
And
$$\align&\sum_{k=0}^\infty\f{\bi{2k}k}{(2k+1)16^k}\(8\sum_{j=0}^k\f{(-1)^j}{(2j+1)^4}+\f{(-1)^k}{(2k+1)^4}\)
\\&\qquad=\f{11}{120}\pi^2\zeta(3)+\f 83\pi \beta(4),\tag1.52
\\&\sum_{k=0}^\infty\f{\bi{2k}k}{(2k+1)16^k}\(\sum_{j=0}^k\f{33}{(2j+1)^5}+\f4{(2k+1)^5}\)
\\&\qquad=\f{35}{288}\pi^3\zeta(3)+\f{1003}{96}\pi\zeta(5),\tag1.53
\\&\sum_{k=0}^\infty\f{\bi{2k}k}{(2k+1)^2(-16)^k}\(110\sum_{j=0}^k\f{(-1)^j}{(2j+1)^4}+29\f{(-1)^k}{(2k+1)^4}\)
\\&\qquad=\f{91}{96}\pi^3\zeta(3)+11\pi^2\beta(4)-\f{301}{192}\pi\zeta(5),\tag1.54
\\&\sum_{k=0}^\infty\f{\bi{2k}k}{(2k+1)^316^k}\(72\sum_{j=0}^k\f{(-1)^j}{(2j+1)^2}-\f{(-1)^k}{(2k+1)^2}\)
\\&\qquad=\f{7}{3}\pi^3G+\f{17}{40}\pi^2\zeta(3),\tag1.55
\endalign$$
$$\align&\sum_{k=0}^\infty\f{\bi{2k}k}{(2k+1)^316^k}\(\sum_{j=0}^k\f{33}{(2j+1)^3}+\f8{(2k+1)^3}\)
\\&\qquad=\f{245}{216}\pi^3\zeta(3)-\f{49}{144}\pi\zeta(5).\tag1.56
\endalign$$
\endproclaim

\Remark. (a) I announced (1.1)-(1.6) first by several messages to {\tt Number Theory Mailing List}
during March-April in 2010. My conjectural identity (1.2)
was confirmed in [G] via applying a Barnes-integrals strategy of the WZ-method. In 2012
Kh. Hessami Pilehrood and T. Hessami Pilehrood [HP] proved my conjectural identity
(1.4) by means of the Hurwitz zeta function. (1.1), (1.3) and (1.5) were recently confirmed by J. Guillera and M. Rogers [GR].
(1.9) was discovered on August 12, 2014. It is known that $\sum_{k=1}^\infty(-1)^{k-1}/(k^3\bi{2k}k)=\f25\zeta(3)$.
A combination of (1.10) and (1.11) yields $\sum_{k=1}^\infty(3H_k-1/k)/(k^2\bi{2k}k)=\zeta(3)$
for which {\tt Mathematica 9} could yield a ``proof" after running the {\tt FullSimplify} command half an hour
(see {\tt http://math.nju.edu.cn/$\sim$zwsun/zeta(3).txt} for my detailed report).
Combining (1.10)-(1.12) we find exact values of
$$\sum_{k=1}^\infty \f1{k^3\bi{2k}k},\ \sum_{k=1}^\infty\f{H_k}{k^2\bi{2k}k}\ \t{and}\ \sum_{k=1}^\infty\f{H_{2k}}{k^2\bi{2k}k}.$$
Note that S. Ramanujan (cf. [BJ]) discovered that
$$\sum_{k=0}^\infty\f{\bi{2k}k}{(2k+1)^28^k}=\f{\pi}{4\sqrt2}\log2+\f{G}{\sqrt2}\ \t{and}\ \sum_{k=0}^\infty\f{\bi{2k}k}{(2k+1)^216^k}=\f{3\sqrt3}4K.$$
In 1985 I.J. Zucker [Z] proved the following remarkable identities:
$$\gather\sum_{k=1}^\infty\f1{k^3\bi{2k}k}=\f{\sqrt3}2\pi K-\f 43\zeta(3),
\ \sum_{k=1}^\infty\f1{k^5\bi{2k}k}=\f{9\sqrt3}8\pi L+\f{\pi^2}9\zeta(3)-\f{19}3\zeta(5),
\\\sum_{k=0}^\infty\f{\bi{2k}k}{(2k+1)^316^k}=\f 7{216}\pi^3,
\ \sum_{k=0}^\infty\f{\bi{2k}k}{(2k+1)^416^k}=\f{\pi\zeta(3)}{12}+\f{27\sqrt3}{64}L.
\endgather$$
Also, (1.19) could be yielded by {\tt Mathematica 9} but it lacks a readable human proof. Concerning (1.20) and (1.21), we remark that (cf. [S14e])
$$\sum_{k=1}^\infty \f{L_{2k}}{k^2\bi{2k}k}=\f{\pi^2}5\ \ \t{and}\ \ \sum_{k=1}^\infty\f{v_k}{k^2\bi{2k}k}=\f25\pi^2.$$
I would like to offer 500 Chinese dollars (RMB) as the prize for the first correct proof of the formula (1.24) for $\zeta(5)$.

(b) L. van Hamme [vH] investigated corresponding $p$-adic congruences for certain hypergeometric series
involving the Gamma function or $\pi=\Gamma(1/2)^2$.
Almost all of my conjectural series were motivated by their $p$-adic analogues that I found first.
For example, (1.9) was motivated by my conjectural congruences
$$\sum_{k=1}^{p-1}\f{\bi{4k}{2k+1}\bi{2k}k}{48^k}\eq\f 5{12}p^2B_{p-2}\l(\f13\r)\pmod{p^2}$$
and $$p^2\sum_{k=1}^{p-1}\f{48^k}{k(2k-1)\bi{4k}{2k}\bi{2k}k}\eq 4\l(\f p3\r)+4p\pmod{p^2}$$
for any prime $p>3$, where $B_{p-2}(x)$ denotes the Bernoulli polynomial of degree $p-2$. Also, (1.31) and (1.43) are related to my conjectural congruences
$$\gather\sum_{k=1}^{p-1}\f{H_k^{(3)}}{k^2\bi{2k}k}\eq\f{29}{45}B_{p-5}\pmod p,\ \sum_{k=1}^{p-1}\f{\bi{2k}k}kH_{k-1}^{(3)}\eq\f2{45}pB_{p-5}\pmod{p^2},
\\\sum_{k=0}^{(p-3)/2}\f{\bi{2k}k}{(2k+1)16^k}\sum_{j=0}^k\f1{(2j+1)^3}\eq\f7{180}\l(\f{-1}p\r)pB_{p-5}\pmod{p^2},
\endgather$$
where $p$ is any prime greater than $3$. 
The reader may consult [S], [S11], [S13b], [S14a], [S14b], [S14c] and [S14e] for many other congruences related to my conjectural series.

\heading{2. Various series for $1/\pi$}\endheading

\proclaim{Conjecture 2} {\rm (i) ([S13b])} Set
$$a_n(x)=\sum_{k=0}^n\bi nk^2\bi{n+k}{k}x^{n-k}\quad(n=0,1,2,\ldots).$$ Then we have
$$\align\sum_{k=0}^\infty\f{13k+4}{96^k}\bi{2k}ka_k(-8)=&\f{9\sqrt2}{2\pi},\tag2.1
\\\sum_{k=0}^\infty\f{290k+61}{1152^k}\bi{2k}ka_k(-32)=&\f{99\sqrt2}{\pi},\tag2.2
\\\sum_{k=0}^\infty\f{962k+137}{3840^k}\bi{2k}ka_k(64)=&\f{252\sqrt5}{\pi}.\tag2.3
\endalign$$

{\rm (ii) (Z. W. Sun [S14c])} For $n=0,1,2,\ldots$ define
$$S_n^{(1)}(x)=\sum_{k=0}^n\bi nk\bi{2k}k\bi{2n-2k}{n-k}x^{n-k},\  S_n^{(2)}(x)=\sum_{k=0}^n\bi{2k}k^2\bi{2n-2k}{n-k}x^{n-k}.$$
Then we have
$$\align\sum_{k=0}^\infty\f{12k+1}{400^k}\bi{2k}kS_k^{(1)}(16)=&\f{25}{\pi},\tag2.4
\\\sum_{k=0}^\infty\f{10k+1}{(-384)^k}\bi{2k}kS_k^{(1)}(-16)=&\f{8\sqrt6}{\pi},\tag2.5
\\\sum_{k=0}^\infty\f{170k+37}{(-3584)^k}\bi{2k}kS_k^{(1)}(64)=&\f{64\sqrt{14}}{3\pi},\tag2.6
\\\sum_{k=0}^\infty\f{476k+103}{3600^k}\bi{2k}kS_k^{(1)}(-64)=&\f{225}{\pi},\tag2.7
\\\sum_{k=0}^\infty\f{140k+19}{4624^k}\bi{2k}kS_k^{(1)}(64)=&\f{289}{3\pi},\tag2.8
\\\sum_{k=0}^\infty\f{1190k+163}{(-4608)^k}\bi{2k}kS_k^{(1)}(-64)=&\f{576\sqrt2}{\pi},\tag2.9
\\\sum_{k=0}^\infty\f{k-1}{72^k}\bi{2k}kS_k^{(2)}(4)=&\f{9}{\pi},\tag$\underline{2.10}$
\\\sum_{k=0}^\infty\f{4k+1}{(-192)^k}\bi{2k}kS_k^{(2)}(4)=&\f{\sqrt{3}}{\pi},\tag2.11
\\\sum_{k=0}^{\infty}\f{k-2}{100^k}\bi{2k}kS_k^{(2)}(6)=&\f{50}{3\pi},\tag2.12
\\\sum_{k=0}^\infty\f{k}{(-192)^k}\bi{2k}kS_k^{(2)}(-8)=&\f{3}{2\pi},\tag2.13
\\\sum_{k=0}^\infty\f{6k-1}{256^k}\bi{2k}kS_k^{(2)}(12)=&\f{8\sqrt3}{\pi},\tag2.14
\\\sum_{k=0}^\infty\f{17k-224}{(-225)^k}\bi{2k}kS_k^{(2)}(-14)=&\f{1800}{\pi},\tag2.15
\\\sum_{k=0}^\infty\f{15k-256}{289^k}\bi{2k}kS_k^{(2)}(18)=&\f{2312}{\pi},\tag2.16
\\\sum_{k=0}^\infty\f{20k-11}{(-576)^k}\bi{2k}kS_k^{(2)}(-32)=&\f{90}{\pi},\tag2.17
\endalign$$
$$\align\sum_{k=0}^\infty\f{10k+1}{(-1536)^k}\bi{2k}kS_k^{(2)}(-32)=&\f{3\sqrt{6}}{\pi},\tag2.18
\\\sum_{k=0}^\infty\f{3k-2}{640^k}\bi{2k}kS_k^{(2)}(36)=&\f{5\sqrt{10}}{\pi},\tag2.19
\\\sum_{k=0}^\infty\f{12k+1}{1600^k}\bi{2k}kS_k^{(2)}(36)=&\f{75}{8\pi},\tag2.20
\\\sum_{k=0}^\infty\f{24k+5}{3136^k}\bi{2k}kS_k^{(2)}(-60)=&\f{49\sqrt{3}}{8\pi},\tag2.21
\\\sum_{k=0}^\infty\f{14k+3}{(-3072)^k}\bi{2k}kS_k^{(2)}(64)=&\f{6}{\pi},\tag2.22
\\\sum_{k=0}^\infty\f{20k-67}{(-3136)^k}\bi{2k}kS_k^{(2)}(-192)=&\f{490}{\pi},\tag2.23
\\\sum_{k=0}^\infty\f{7k-24}{3200^k}\bi{2k}kS_k^{(2)}(196)=&\f{125\sqrt2}{\pi},\tag2.24
\\\sum_{k=0}^\infty\f{5k-32}{(-6336)^k}\bi{2k}kS_k^{(2)}(-392)=&\f{495}{2\pi},\tag2.25
\\\sum_{k=0}^\infty\f{66k-427}{6400^k}\bi{2k}kS_k^{(2)}(396)=&\f{1000\sqrt{11}}{\pi},\tag2.26
\\\sum_{k=0}^\infty\f{34k-7}{(-18432)^k}\bi{2k}kS_k^{(2)}(-896)=&\f{54\sqrt{2}}{\pi},\tag2.27
\\\sum_{k=0}^\infty\f{24k-5}{18496^k}\bi{2k}kS_k^{(2)}(900)=&\f{867}{16\pi}.\tag2.28
\endalign$$
\endproclaim
\Remark. (i) Those $a_n(1)\ (n=0,1,2,\ldots)$ were first introduced by R. Ap\'ery in his study of the irrationality of $\zeta(2)$ and $\zeta(3)$.
Identities related to the form $\sum_{k=0}^\infty (bk+c)\bi{2k}ka_k(1)/m^k=C/\pi$
were first studied by T. Sato in 2002.

(ii) I introduced the polynomials $S_n^{(1)}(x)$ and $S_n^{(2)}(x)$ during March 27-28, 2011.
(2.4)-(2.23) and (2.24)-(2.28) were discovered during March 27-31, 2011 and  Jan. 23-24, 2012 respectively.
By {\tt Mathematica}, we have
$$S_n^{(1)}(-1)=\cases\bi{n}{n/2}^2&\t{if}\ 2\mid n,\\0&\t{if}\ 2\nmid n.\endcases$$
I also noted that
$$S_n^{(1)}(1)=\sum_{k=0}^{\lfloor n/2\rfloor}\bi n{2k}\bi{2k}k^24^{n-2k}.$$
Identities of the form
$$\sum_{n=0}^\infty \f{bn+c}{m^n}\bi{2n}n\sum_{k=0}^{\lfloor n/2\rfloor}\bi n{2k}\bi{2k}k^24^{n-2k}=\f C{\pi}$$
were recently investigated in [CC].

(iii) In [S14c] I proved the following three identities via Ramanujan-type series for $1/\pi$ (cf. [B, pp.\,353-354]).
$$\align \sum_{k=0}^\infty\f{k}{128^k}\bi{2k}kS_k^{(2)}(4)=&\f{\sqrt2}{\pi},
\\\sum_{k=0}^\infty\f{8k+1}{576^k}\bi{2k}kS_k^{(2)}(4)=&\f{9}{2\pi},
\\\sum_{k=0}^\infty\f{8k+1}{(-4032)^k}\bi{2k}kS_k^{(2)}(4)=&\f{9\sqrt7}{8\pi}.
\endalign$$
In December 2011 A. Meurman [M] confirmed my conjectural (2.10).

\proclaim{Conjecture 3} {\rm {(i) (Discovered on April 1, 2011)} Set
$$W_n(x):=\sum_{k=0}^n\bi{n+k}{2k}\bi{2k}k^2\bi{2n-2k}{n-k}x^{-(n+k)}\ \ (n=0,1,2,\ldots).$$
Then
$$\align\sum_{k=0}^\infty (8k+3)W_k(-8)=&\f{28\sqrt3}{9\pi},\tag3.1
\\\sum_{k=0}^\infty(8k+1)W_k(12)=&\f{26\sqrt3}{3\pi},\tag3.2
\\\sum_{k=0}^\infty (24k+7)W_k(-16)=&\f{8\sqrt3}{\pi},\tag3.3
\\\sum_{k=0}^\infty(360k+51)W_k(20)=&\f{210\sqrt3}{\pi},\tag3.4
\\\sum_{k=0}^\infty(21k+5)W_k(-28)=&\f{63\sqrt2}{8\pi},\tag3.5
\\\sum_{k=0}^\infty(7k+1)W_k(32)=&\f{11\sqrt2}{3\pi},\tag3.6
\endalign$$
$$\align\sum_{k=0}^\infty(195k+31)W_k(-100)=&\f{275\sqrt6}{8\pi},\tag3.7
\\\sum_{k=0}^\infty(39k+5)W_k(104)=&\f{91\sqrt6}{12\pi},\tag3.8
\\\sum_{k=0}^\infty(2856k+383)W_k(-196)=&\f{637\sqrt3}{\pi},\tag3.9
\\\sum_{k=0}^\infty(14280k+1681)W_k(200)=&\f{3350\sqrt3}{\pi}.\tag3.10
\endalign$$

 {\rm (ii) (Discovered during April 7-10, 2011 and Oct. 6-7, 2012; (3.18), (3.24)-(3.25) and (3.28) appeared in [S13b])} Define
$$f_n^+(x):=\sum_{k=0}^n\bi nk^2\bi{2k}nx^{2k-n}$$
and
$$f_n^-(x):=\sum_{k=0}^n\bi nk^2\bi{2k}n(-1)^kx^{2k-n}$$
for $n=0,1,2,\ldots$. Then
$$\align \sum_{k=0}^\infty\f{19k+3}{240^k}\bi{2k}kf^+_k(6)=&\f{35\sqrt6}{4\pi},\tag3.11
\\\sum_{k=0}^\infty\f{135k+8}{289^k}\bi{2k}kf^+_k(14)=&\f{6647}{14\pi},\tag3.12
\\\sum_{k=0}^\infty\f{297k+41}{2800^k}\bi{2k}kf^+_k(14)=&\f{325\sqrt{14}}{8\pi},\tag3.13
\\\sum_{k=0}^\infty\f{770k+79}{576^k}\bi{2k}kf^+_k(21)=&\f{468\sqrt7}{\pi},\tag3.14
\\\sum_{k=0}^\infty\f{209627k+22921}{46800^k}\bi{2k}kf^+_k(36)=&\f{58275\sqrt{26}}{4\pi},\tag3.15
\endalign$$
$$\align\sum_{k=0}^\infty\f{322k+41}{2304^k}\bi{2k}kf^+_k(45)=&\f{3456\sqrt{7}}{35\pi},\tag3.16
\\\sum_{k=0}^\infty\f{205868k+18903}{439280^k}\bi{2k}kf^+_k(76)=&\f{1112650\sqrt{19}}{81\pi},\tag3.17
\\\sum_{k=0}^\infty\f{8851815k+1356374}{(-29584)^k}\bi{2k}kf^+_k(175)=&\f{1349770\sqrt{7}}{\pi},\tag3.18
\\\sum_{k=0}^\infty\f{12980k-2303}{5616^k}\bi{2k}kf^+_k(300)=&\f{34398\sqrt{3}}{\pi},\tag3.19
\\\sum_{k=0}^\infty\f{1391k+21}{28880^k}\bi{2k}kf^+_k(1156)=&\f{229957\sqrt{10}}{324\pi},\tag3.20
\\\sum_{k=0}^\infty\f{68572k-34329}{20400^k}\bi{2k}kf^+_k(1176)=&\f{82450\sqrt{51}}{\pi},\tag3.21
\\\sum_{k=0}^\infty\f{930886k-159493}{243360^k}\bi{2k}kf^+_k(12321)=&\f{5636826\sqrt{95}}{19\pi},\tag3.22
\\\sum_{k=0}^\infty\f{182k+51}{48^k}\bi{2k}kf_k^-\l(\f{15}{16}\r)=&\f{552}{5\pi},\tag3.23
\\\sum_{k=0}^\infty\f{1054k+233}{480^k}\bi{2k}kf^-_k(8)=&\f{520}{\pi},\tag$\underline{3.24}$
\\\sum_{k=0}^\infty\f{224434k+32849}{5760^k}\bi{2k}kf^-_k(18)=&\f{93600}{\pi},\tag3.25
\endalign$$
$$\align
\\\sum_{k=0}^\infty\f{170k+41}{(-48)^k}\bi{2k}kf^-_k\l(\f 98\r)=&\f{78\sqrt6}{\pi},\tag3.26
\\\sum_{k=0}^\infty\f{15470k+1063}{(-288)^k}\bi{2k}kf^-_k\l(\f {225}{16}\r)=&\f{37044\sqrt2}{\pi}.\tag3.27
\endalign$$

{\rm (iii) ([S13b])}}  Define
$g_n(x)=\sum_{k=0}^n\bi nk^2\bi{2k}kx^k$ for $n=0,1,2,\ldots.$
Then
$$\sum_{k=0}^\infty\f{16k+5}{18^{2k}}\bi{2k}kg_k(-20)=\f{189}{25\pi}.\tag3.28$$
We also have
$$\sum_{n=0}^\infty\f{21n+1}{64^n}\sum_{k=0}^n\bi nk\bi{2k}n\bi{2k}k\bi{2n-2k}{n-k}3^{2k-n}=\f{64}{\pi}.\tag3.29$$
\endproclaim
\Remark. (a)  As May 20 is the day for Nanjing University, I offered \$520 (520 US dollars) for the first correct proof of (3.24).
Later, M. Rogers and A. Straub [RS] won the prize, and they also discussed other series in Conjecture 3(ii).

(b) For $n=0,1,2,\ldots$ define
$$f_n(x):=\sum_{k=0}^n\bi{n}k^2\bi{2k}nx^k=\sum_{k=0}^n\bi nk\bi{2k}k\bi{k}{n-k}x^k.$$
Then
$$f_n(1)=\sum_{k=0}^n\bi nk^3,\ \ \ f_n^+(x)=x^{-n}f_n(x^2)\ \ \ \t{and}\ \ \ f_n^-(x)=x^{-n}f_n(-x^2).$$
By [S14d], $\sum_{k=0}^n\bi nk(-1)^k((-1)^kf_k(x))=g_n(x)$. Thus, by
the technique in Section 5, each of (3.11)-(3.27) has an equivalent
form in term of $g_k(x)$. Below are equivalent versions of
(3.11)-(3.15), (3.17)-(3.18) and (3.24)-(3.25):
$$\align
\sum_{k=0}^\infty\f{720k+113}{38^{2k}}\bi{2k}kg_k(36)=&\f{2527\sqrt{15}}{12\pi},\tag3.11$'$
\\\sum_{k=0}^\infty\f{17k+1}{4050^k}\bi{2k}kg_k(196)=&\f{15525}{98\sqrt7\,\pi},\tag3.12$'$
\\\sum_{k=0}^\infty\f{3920k+541}{198^{2k}}\bi{2k}kg_k(196)=&\f{42471}{8\sqrt7\,\pi},\tag3.13$'$
\\\sum_{k=0}^\infty\f{2352k+241}{110^{2k}}\bi{2k}kg_k(441)=&\f{39325}{6\sqrt3\,\pi},\tag3.14$'$
\\\sum_{k=0}^\infty\f{18139680k+1983409}{1298^{2k}}\bi{2k}kg_k(1296)=&\f{109091059}{12\sqrt2\,\pi},\tag3.15$'$
\endalign$$
$$\align\sum_{k=0}^\infty\f{944607040k+86734691}{5778^{2k}}\bi{2k}kg_k(5776)=&\f{1071111195\sqrt{95}}{38\pi},\tag3.17$'$
\\\sum_{k=0}^\infty\f{35819000k+5488597}{(-5177196)^k}\bi{2k}kg_k(30625)=&\f{3315222\sqrt{19}}{\pi},\tag3.18$'$
\\\sum_{k=0}^\infty\f{5440k+1201}{62^{2k}}\bi{2k}kg_k(-64)=&\f{12493\sqrt{15}}{18\pi},\tag3.24$'$
\\\sum_{k=0}^\infty\f{1505520k+220333}{322^{2k}}\bi{2k}kg_k(-324)=&\f{1684865\sqrt{5}}{6\pi}.\tag3.25$'$
\endalign$$
Note that [CTYZ] contains some series for $1/\pi$ involving $f_k=f_k(1)$ or $g_k=g_k(1)$.

(c) Observe that
$$\sum_{k=0}^n\bi nk\bi{2k}n\bi{2k}k\bi{2n-2k}{n-k}(-1)^{n-k}
=\sum_{k=0}^n\bi nk^2\bi{2k}k\bi{2n-2k}{n-k},$$
which can be proved by obtaining the same recurrence relation for both sides via the Zeilberger algorithm.

\medskip

\proclaim{Conjecture 4 {\rm (Discovered during April 23-25 and May 7-16, 2011)}}

{\rm (i)} We have
$$\align\sum_{n=0}^\infty\f{8n+1}{9^n}\sum_{k=0}^n\bi{-1/3}k^2\bi{-2/3}{n-k}^2=&\f {3\sqrt3}{\pi},\tag$\underline{4.1}$
\\\sum_{n=0}^\infty\f{(2n-1)(-3)^n}{16^n}\sum_{k=0}^n\bi{2k}k\bi{2(n-k)}{n-k}\bi{-1/3}k\bi{-2/3}{n-k}=&\f {16}{\sqrt3\,\pi},\tag$\underline{4.2}$
\\\sum_{n=0}^\infty\f{10n+3}{16^n}\sum_{k=0}^n\bi{2k}k\bi{2(n-k)}{n-k}\bi{-1/3}k\bi{-2/3}{n-k}=&\f {16\sqrt3}{5\pi},\tag$\underline{4.3}$
\\\sum_{n=0}^\infty\f{8n+1}{(-20)^n}\sum_{k=0}^n\bi{2k}k\bi{2(n-k)}{n-k}\bi{-1/3}k\bi{-2/3}{n-k}=&\f {4\sqrt3}{\pi},\tag$\underline{4.4}$
\\\sum_{n=0}^\infty\f{168n+29}{108^n}\sum_{k=0}^n\bi{2k}k\bi{2(n-k)}{n-k}\bi{-1/3}k\bi{-2/3}{n-k}=&\f {324\sqrt3}{7\pi},\tag$\underline{4.5}$
\\\sum_{n=0}^\infty\f{162n+23}{(-112)^n}\sum_{k=0}^n\bi{2k}k\bi{2(n-k)}{n-k}\bi{-1/3}k\bi{-2/3}{n-k}=&\f {48\sqrt3}{\pi}.\tag$\underline{4.6}$
\endalign$$
Also,
$$\align\sum_{n=0}^\infty\f{(n-2)(-2)^n}{9^n}\sum_{k=0}^n\bi{2k}k\bi{2(n-k)}{n-k}\bi{-1/4}k\bi{-3/4}{n-k}
=&\f{6\sqrt3}{\pi},\tag$\underline{4.7}$
\\\sum_{n=0}^\infty\f{16n+5}{12^n}\sum_{k=0}^n\bi{2k}k\bi{2(n-k)}{n-k}\bi{-1/4}k\bi{-3/4}{n-k}=&\f 8{\pi},\tag$\underline{4.8}$
\\\sum_{n=0}^\infty\f{12n+1}{(-16)^n}\sum_{k=0}^n\bi{2k}k\bi{2(n-k)}{n-k}\bi{-1/4}k\bi{-3/4}{n-k}=&\f {32}{3\pi},\tag$\underline{4.9}$
\endalign$$
and
$$\align\sum_{n=0}^\infty\f{(81n+32)8^n}{49^n}\sum_{k=0}^n\bi{2k}k\bi{2(n-k)}{n-k}\bi{-1/4}k\bi{-3/4}{n-k}=&\f {14\sqrt7}{\pi},\tag$\underline{4.10}$
\\\sum_{n=0}^\infty\f{n(-8)^n}{81^n}\sum_{k=0}^n\bi{2k}k\bi{2(n-k)}{n-k}\bi{-1/4}k\bi{-3/4}{n-k}=&\f {5}{4\pi}\,?\tag$\underline{4.11}$
\\\sum_{n=0}^\infty\f{324n+43}{320^n}\sum_{k=0}^n\bi{2k}k\bi{2(n-k)}{n-k}\bi{-1/4}k\bi{-3/4}{n-k}=&\f {128}{\pi},\tag$\underline{4.12}$
\\\sum_{n=0}^\infty\f{320n+39}{(-324)^n}\sum_{k=0}^n\bi{2k}k\bi{2(n-k)}{n-k}\bi{-1/4}k\bi{-3/4}{n-k}=&\f {648}{5\pi}.\tag$\underline{4.13}$
\endalign$$

{\rm (ii)} We have
$$\sum_{n=0}^\infty\f{3n-1}{2^n}\sum_{k=0}^n\bi {-1/3}k\bi{-2/3}{n-k}\bi{-1/6}k\bi{-5/6}{n-k}=\f{3\sqrt6}{2\pi}.\tag4.14$$
If we set
$$\align a_n=&\sum_{k=0}^n(-1)^k\bi{-1/3}k^2\bi{-2/3}{n-k}
\\=&\sum_{k=0}^n(-1)^k\bi{-2/3}k^2\bi{-1/3}{n-k}
\\=&\f{(-4)^n}{\bi{2n}n}\sum_{k=0}^n\bi{-2/3}k\bi{-1/3}{n-k}\bi{-1/6}k\bi{-5/6}{n-k},
\endalign$$
then
$$\align\sum_{n=0}^\infty\f{3n-2}{(-5)^n}\bi{2n}na_n=&\f{3\sqrt{15}}{\pi},\tag$\underline{4.15}$
\\\sum_{n=0}^\infty\f{32n+1}{(-100)^n}9^n\bi{2n}na_n=&\f{50}{\sqrt3\,\pi},\tag$\underline{4.16}$
\\\sum_{n=0}^\infty\f{81n+13}{50^n}\bi{2n}na_n=&\f{75\sqrt3}{4\pi},\tag$\underline{4.17}$
\endalign$$
$$\align\sum_{n=0}^\infty\f{96n+11}{(-68)^n}\bi{2n}na_n=&\f{6\sqrt{51}}{\pi},\tag$\underline{4.18}$
\\\sum_{n=0}^\infty\f{15n+2}{121^n}\bi{2n}na_n=&\f{363\sqrt{15}}{250\pi},\tag$\underline{4.19}$
\\\sum_{n=0}^\infty\f{160n+17}{(-324)^n}\bi{2n}na_n=&\f{16}{\sqrt3\,\pi},\tag$\underline{4.20}$
\\\sum_{n=0}^\infty\f{6144n+527}{(-4100)^n}\bi{2n}na_n=&\f{150\sqrt{123}}{\pi},\tag$\underline{4.21}$
\\\sum_{n=0}^\infty\f{1500000n+87659}{(-1000004)^n}\bi{2n}na_n=&\f{16854\sqrt{267}}{\pi}.\tag$\underline{4.22}$
\endalign$$

{\rm (iii)} For $n=0,1,2,\ldots$ set
$$\align b_n
=&\sum_{k=0}^n(-1)^k\bi{-1/4}k^2\bi{-3/4}{n-k}
\\=&\sum_{k=0}^n(-1)^k\bi{-3/4}k^2\bi{-1/4}{n-k}
\\=&\f{(-4)^n}{\bi{2n}n}\sum_{k=0}^n\bi{-1/8}k\bi{-5/8}k\bi{-3/8}{n-k}\bi{-7/8}{n-k}.
\endalign$$
Then
$$\align\sum_{n=0}^\infty\f{16n+1}{(-20)^n}\bi{2n}nb_n=&\f{4\sqrt5}{\pi},\tag$\underline{4.23}$
\\\sum_{n=0}^\infty\f{(3n-1)4^n}{(-25)^n}\bi{2n}nb_n=&\f{25}{3\sqrt3\,\pi},\tag$\underline{4.24}$
\\\sum_{n=0}^\infty\f{6n+1}{32^n}\bi{2n}nb_n=&\f{8\sqrt{6}}{9\pi},\tag$\underline{4.25}$
\\\sum_{n=0}^\infty\f{81n+23}{49^n}8^n\bi{2n}nb_n=&\f{49}{2\pi},\tag$\underline{4.26}$
\\\sum_{n=0}^\infty\f{192n+19}{(-196)^n}\bi{2n}nb_n=&\f{196}{3\pi},\tag$\underline{4.27}$
\endalign$$
and
$$\align
\sum_{n=0}^\infty\f{162n+17}{320^n}\bi{2n}nb_n=&\f{16\sqrt{10}}{\pi},\tag$\underline{4.28}$
\\\sum_{n=0}^\infty\f{1296n+113}{(-1300)^n}\bi{2n}nb_n=&\f{100\sqrt{13}}{\pi},\tag$\underline{4.29}$
\\\sum_{n=0}^\infty\f{4802n+361}{9600^n}\bi{2n}nb_n=&\f{800\sqrt{2}}{\pi},\tag$\underline{4.30}$
\\\sum_{n=0}^\infty\f{162n+11}{39200^n}\bi{2n}nb_n=&\f{19600}{121\sqrt{22}\,\pi}.\tag$\underline{4.31}$
\endalign$$

{\rm (iv)} For $n=0,1,2,\ldots$ set
$$\align c_n:
=&\sum_{k=0}^n(-1)^k\bi{-1/6}k^2\bi{-5/6}{n-k}
\\=&\sum_{k=0}^n(-1)^k\bi{-5/6}k^2\bi{-1/6}{n-k}
\\=&\f{(-4)^n}{\bi{2n}n}\sum_{k=0}^n\bi{-1/12}k\bi{-7/12}k\bi{-5/12}{n-k}\bi{-11/12}{n-k}.
\endalign$$
Then we have
$$\align\sum_{n=0}^\infty\f{125n+13}{121^n}\bi{2n}nc_n=&\f{121}{2\sqrt3\,\pi},\tag$\underline{4.32}$
\\\sum_{n=0}^\infty\f{(125n-8)16^n}{(-189)^n}\bi{2n}nc_n=&\f{27\sqrt7}{\pi},\tag$\underline{4.33}$
\\\sum_{n=0}^\infty\f{(125n+24)27^n}{392^n}\bi{2n}nc_n=&\f{49}{\sqrt2\,\pi},\tag$\underline{4.34}$
\\\sum_{n=0}^\infty\f{512n+37}{(-2052)^n}\bi{2n}nc_n=&\f{27\sqrt{19}}{\pi}.\tag$\underline{4.35}$
\\\sum_{n=0}^\infty\f{(512n+39)27^n}{(-2156)^n}\bi{2n}nc_n=&\f{49\sqrt{11}}{\pi},\tag$\underline{4.36}$
\\\sum_{n=0}^\infty\f{(1331n+109)2^n}{1323^n}\bi{2n}nc_n=&\f{1323}{4\pi}.\tag$\underline{4.37}$
\endalign$$
\endproclaim
\Remark. (i) I [S11a] proved the following three identities:
$$\align \sum_{n=0}^\infty\f n{4^n}\sum_{k=0}^n\bi{-1/4}k^2\bi{-3/4}{n-k}^2=&\f{4\sqrt3}{9\pi},
\\\sum_{n=0}^\infty\f{9n+2}{(-8)^n}\sum_{k=0}^n\bi{-1/4}k^2\bi{-3/4}{n-k}^2=&\f 4{\pi},
\\\sum_{n=0}^\infty\f{9n+1}{64^n}\sum_{k=0}^n\bi{-1/4}k^2\bi{-3/4}{n-k}^2=&\f {64}{7\sqrt7\,\pi}.
\endalign$$
On May 15, 2011 I observed that if $x+y+1=0$ then
$$\sum_{k=0}^n(-1)^k\bi xk^2\bi{y}{n-k}=\sum_{k=0}^n(-1)^k\bi yk^2\bi{x}{n-k}$$
which can be easily proved since both sides satisfy the same recurrence relation by the Zeilberger algorithm.
Also,
$$\align\sum_{k=0}^n\bi{-1/3}k\bi{-2/3}{n-k}\bi{-1/6}k\bi{-5/6}{n-k}=&\sum_{k=0}^n\bi nk\bi{2k}k\f{a_k}{4^k},
\\\sum_{k=0}^n\bi{-1/8}k\bi{-3/8}k\bi{-5/8}{n-k}\bi{-7/8}{n-k}=&\sum_{k=0}^n\bi nk\bi{2k}k\f{b_k}{4^k},
\\\sum_{k=0}^n\bi{-1/12}k\bi{-5/12}k\bi{-7/12}{n-k}\bi{-11/12}{n-k}=&\sum_{k=0}^n\bi nk\bi{2k}k\f{c_k}{4^k}.
\endalign$$
Thus, each of (4.15)-(4.37) has an equivalent form since
$$\sum_{n=0}^\infty\f{bn+c}{m^n}\sum_{k=0}^n\bi nk(-1)^kf(k)=
\f m{(m-1)^2}\sum_{k=0}^\infty\f{bmk+b+(m-1)c}{(1-m)^k}f(k)$$
if both series in the equality converge absolutely. For example, (4.22) holds if and only if
$$\sum_{n=0}^\infty\f{16854n+985}{(-250000)^n}\sum_{k=0}^n\bi{-1/3}k\bi{-2/3}{n-k}\bi{-1/6}k\bi{-5/6}{n-k}=\f{4500000}{89\sqrt{267}\,\pi}.$$

{(ii)} (4.1) and (4.14) appeared as conjectures in [S13b].
In December 2011 G. Almkvist and A. Aycock  released the preprint [AA] in which they proved all the conjectured formulas in Conj. 4
except (4.14), with the right-hand side of (4.11) corrected as $162/(49\sqrt7\,\pi)$.

\proclaim{Conjecture 5 {\rm (Z. W. Sun [S14c])}} Define
$$s_n(x):=\sum_{k=0}^n\bi nk\bi{n+2k}{2k}\bi{2k}kx^{-(n+k)}\quad\t{for}\ n=0,1,2,\ldots$$
Then
$$\align\sum_{k=0}^\infty(7k+2)\bi{2k}ks_k(-9)&=\f{9\sqrt3}{5\pi},\tag5.1
\\\sum_{k=0}^\infty(9k+2)\bi{2k}ks_k(-20)&=\f{4}{\pi},\tag5.2
\\\sum_{k=0}^\infty(95k+13)\bi{2k}ks_k(36)&=\f{18\sqrt{15}}{\pi},\tag5.3
\\\sum_{k=0}^\infty(310k+49)\bi{2k}ks_k(-64)&=\f{32\sqrt{15}}{\pi},\tag5.4
\\\sum_{k=0}^\infty(495k+53)\bi{2k}ks_k(196)&=\f{70\sqrt{7}}{\pi},\tag5.5
\\\sum_{k=0}^\infty(13685k+1474)\bi{2k}ks_k(-324)&=\f{1944\sqrt5}{\pi},\tag5.6
\\\sum_{k=0}^\infty(3245k+268)\bi{2k}ks_k(1296)&=\f{1215}{\sqrt2\,\pi},\tag5.7
\\\sum_{k=0}^\infty(6420k+443)\bi{2k}ks_k(5776)&=\f{1292\sqrt{95}}{9\pi}.\tag5.8
\endalign$$
Also,
$$\align\sum_{n=0}^\infty\f{357n+103}{2160^n}\bi{2n}n\sum_{k=0}^n\bi nk\bi{n+2k}{2k}\bi{2k}k(-324)^{n-k}=&\f{90}{\pi},\tag5.9
\\\sum_{n=0}^\infty\f{n}{3645^n}\bi{2n}n\sum_{k=0}^n\bi nk\bi{n+2k}{2k}\bi{2k}k486^{n-k}=&\f{10}{3\pi}.\tag5.10
\endalign$$
\endproclaim
\Remark. (5.1)-(5.9) and (5.10) were discovered during June 16-17, 2011 and on Jan. 18, 2012 respectively. I would like to offer \$90 for the first rigorous proof of (5.9) (which first appeared in Conjecture 1.7 of [S13b]), and \$105 for the first complete proof of my following related conjecture: For any prime $p>5$, we have
$$\align&\sum_{n=0}^{p-1}\f{357n+103}{2160^n}\bi{2n}n\sum_{k=0}^n\bi nk\bi{n+2k}{2k}\bi{2k}k(-324)^{n-k}
\\\quad\qquad&\eq p\l(\f{-1}p\r)\l(54+49\l(\f p{15}\r)\r)\pmod{p^2},
\endalign$$
and
$$\align&\sum_{n=0}^{p-1}\f{\bi{2n}n}{2160^n}\sum_{k=0}^n\bi nk\bi{n+2k}{2k}\bi{2k}k(-324)^{n-k}
\\\eq&\cases 4x^2-2p\,(\mo\ p^2)&\t{if}\, (\f{-1}p)=(\f p3)=(\f p5)=(\f p7)=1,\, p=x^2+105y^2,
\\2x^2-2p\,(\mo\ p^2)&\t{if}\, (\f{-1}p)=(\f p7)=1,\,(\f p3)=(\f p5)=-1,\, 2p=x^2+105y^2,
\\2p-12x^2\,(\mo\ p^2)&\t{if}\,(\f{-1}p)=(\f p3)=(\f p5)=(\f p7)=-1,\, p=3x^2+35y^2,
\\2p-6x^2\,(\mo\ p^2)&\t{if}\, (\f{-1}p)=(\f p7)=-1,\, (\f p3)=(\f p5)=1,\, 2p=3x^2+35y^2,
\\20x^2-2p\,(\mo\ p^2)&\t{if}\, (\f{-1}p)=(\f p5)=1,\,(\f p3)=(\f p7)=-1,\, p=5x^2+21y^2,
\\10x^2-2p\,(\mo\ p^2)&\t{if}\, (\f{-1}p)=(\f p3)=1,\, (\f p5)=(\f p7)=-1,\, 2p=5x^2+21y^2,
\\28x^2-2p\,(\mo\ p^2)&\t{if}\, (\f{-1}p)=(\f p5)=-1,\,(\f p3)=(\f p7)=1,\, p=7x^2+15y^2,
\\14x^2-2p\,(\mo\ p^2)&\t{if}\, (\f{-1}p)=(\f p3)=-1,\, (\f p5)=(\f p7)=1,\, 2p=7x^2+15y^2,
\\0\pmod{p^2}&\t{if}\ (\f{-105}p)=-1.
\endcases
\endalign$$
(Note that the imaginary quadratic field $\Q(\sqrt{-105})$ has class number eight.)
In fact, for all series for $1/\pi$ that I found, I had such conjectures on congruences.
See [S13b] for my philosophy about series for $1/\pi$.
\medskip

\proclaim{Conjecture 6} {\rm (i) ([S13b])} We have
$$\align\sum_{n=0}^\infty\f{114n+31}{26^{2n}}\bi{2n}n\sum_{k=0}^n\bi nk^2\bi{n+k}k(-27)^k=&\f{338\sqrt3}{11\pi},\tag$\underline{6.1}$
\\\sum_{n=0}^\infty\f{930n+143}{28^{2n}}\bi{2n}n\sum_{k=0}^n\bi nk^2\bi{n+k}k27^k=&\f{980\sqrt3}{\pi}.\tag$\underline{6.2}$
\endalign$$

{\rm (ii) ((6.3)-(6.7) and (6.8)-(6.13) were discovered on Jan. 12, 2012 and Nov. 16, 2014 respectively)} Set
$$P_n(x)=\sum_{k=0}^n\bi{2k}k^2\bi{k}{n-k}x^{n-k}\ \quad\t{for}\ \ n=0,1,2,\ldots.$$
Then we have
$$\align\sum_{k=0}^\infty\f{14k+3}{8^{2k}}\bi{2k}kP_k(-7)=&\f{16\sqrt7}{3\pi},\tag6.3
\\\sum_{k=0}^\infty\f{255k+56}{13^{2k}}\bi{2k}kP_k(14)=&\f{2028}{7\pi},\tag6.4
\\\sum_{k=0}^\infty\f{308k+59}{20^{2k}}\bi{2k}kP_k(21)=&\f{250\sqrt7}{3\pi},\tag6.5
\\\sum_{k=0}^\infty\f{1932k+295}{44^{2k}}\bi{2k}kP_k(45)=&\f{363\sqrt7}{\pi},\tag6.6
\\\sum_{k=0}^\infty\f{890358k+97579}{176^{2k}}\bi{2k}kP_k(-175)=&\f{116160\sqrt7}{\pi},\tag6.7
\\\sum_{k=0}^\infty\f{130k+41}{384^k}\bi{2k}kP_k(-196)=&\f{112}{\pi},\tag6.8
\\\sum_{k=0}^\infty\f{46k+13}{(-400)^k}\bi{2k}kP_k(196)=&\f{175\sqrt6}{9\pi},\tag6.9
\\\sum_{k=0}^\infty\f{510k+143}{784^k}\bi{2k}kP_k(-396)=&\f{294\sqrt2}{\pi},\tag6.10
\\\sum_{k=0}^\infty\f{42k+11}{(-800)^k}\bi{2k}kP_k(396)=&\f{75}{2\pi},\tag6.11
\\\sum_{k=0}^\infty\f{1848054k+309217}{78400^k}\bi{2k}kP_k(-39204)=&\f{970200}{\pi},\tag6.12
\\\sum_{k=0}^\infty\f{171465k+28643}{(-78416)^k}\bi{2k}kP_k(39204)=&\f{16731\sqrt{29}}{\pi}.\tag6.13
\endalign$$

{\rm (iii) (Z. W. Sun [S14b, (1.8)])} Define
$$s_n=\sum_{k=0}^n5^k\bi{2k}k^2\bi{2(n-k)}{n-k}^2\bigg/\bi nk\quad\t{for}\ \ n=0,1,2,\ldots.$$
Then we have
$$\sum_{k=0}^\infty\f{28k+5}{576^k}\bi{2k}ks_k=\f{9}{\pi}(2+\sqrt2).\tag6.14$$
\endproclaim
\Remark. (i) W. Zudilin [Zu] confirmed (6.1) and (6.2).

(ii) (6.14) was discovered on Jan. 14, 2012. It is known that $\bi nk\mid\bi{2k}k\bi{2(n-k)}{n-k}$ for all $k=0,\ldots,n$.
Recall that the Catalan-Larcombe-French numbers $P_0,P_1,\ldots$ are given by
$$P_n=\sum_{k=0}^n\f{\bi{2k}k^2\bi{2(n-k)}{n-k}^2}{\bi nk}=2^nP_n(-4)=2^n\sum_{k=0}^{\lfloor n/2\rfloor}\bi n{2k}\bi{2k}k^24^{n-2k}$$
and these numbers satisfy the recurrence relation
$$(k+1)^2P_{k+1}=(24k(k+1)+8)P_k-128k^2P_{k-1}\ (k=1,2,3,\ldots).$$
Note  that
$$\sum_{k=0}^n(-1)^k\f{\bi{2k}k^2\bi{2(n-k)}{n-k}^2}{\bi nk}=\cases4^n\bi n{n/2}^2&\t{if}\ 2\mid n,
\\0&\t{if}\ 2\nmid n.\endcases$$
The sequence $\{s_n\}_{n\gs0}$ can also be defined by
$s_0=1,\ s_1=24,\ s_2=976$ and the recurrence relation
$$\align &51200(n+1)^2(n+3)s_n-1920(4n^3+24n^2+46n+29)s_{n+1}
\\&+8(n+2)(41n^2+205n+255)s_{n+2}-3(n+2)(n+3)^2s_{n+3}=0.
\endalign$$
\medskip

A sequence of polynomials $\{P_n(q)\}_{n\gs0}$ with integer coefficients is said to be {\it $q$-logconvex} if
for each  $n=1,2,3,\ldots$ all the coefficients of the polynomial $P_{n-1}(q)P_{n+1}(q)-P_n(q)^2\in\Z[q]$ are nonnegative.
In view of Conjectures 2 and 3, on May 7, 2011 I conjectured that $\{P_n(q)\}_{n\gs0}$ is $q$-logconvex if $P_n(q)$ has one of the following forms:
$$\gather\sum_{k=0}^n\bi nk^2\bi{n+k}kq^k,\ \ \sum_{k=0}^n\bi nk\bi{2k}k\bi{2(n-k)}{n-k}q^k,
\\\sum_{k=0}^n\bi nk^2\bi{2k}k\bi{2(n-k)}{n-k}q^k,
\ \ \sum_{k=0}^n\bi{n+k}{2k}\bi{2k}k^2\bi{2(n-k)}{n-k}q^k.
\endgather$$
For polynomials of the third form, this was later confirmed by D.Q.J. Dou and A.X.Y. Ren [DR].

\heading{3. Series for $1/\pi$ involving generalized central trinomial coefficients}\endheading

For $b,c\in\Z$, the {\it generalized central trinomial coefficient}
$T_n(b,c)$ denotes the coefficient of $x^n$ in
the expansion of $(x^2+bx+c)^n$. It is easy to see that
$$T_n(b,c)=\sum_{k=0}^{\lfloor n/2\rfloor}\bi n{2k}\bi{2k}kb^{n-2k}c^k
=\sum_{k=0}^{\lfloor n/2\rfloor}\bi nk\bi{n-k}k b^{n-2k}c^k.$$ An
efficient way to compute $T_n(b,c)$ is to use the initial values
$$T_0(b,c)=1,\ \ T_1(b,c)=b,$$ and the recursion
$$(n+1)T_{n+1}(b,c)=(2n+1)bT_n(b,c)-n(b^2-4c)T_{n-1}(b,c)\ \ (n=1,2,\ldots).$$
In view of the Laplace-Heine asymptotic formula for Legendre
polynomials, I [S14a] noted that for any positive reals $b$
and $c$ we have
$$T_n(b,c)\sim\f{(b+2\sqrt{c})^{n+1/2}}{2\root{4}\of c\sqrt{n\pi}}$$
as $n\to+\infty$.
For any real number $b$ and $c<0$, I [S14a] conjectured that $\lim_{n\to\infty}\root n\of{|T_n(b,c)|}=\sqrt{b^2-4c}$,
which was later confirmed by S. Wagner [Wa].

In Jan.-Feb. 2011, I introduced a number of series for $1/\pi$ of
the following new types with $a,b,c,d,m$ integers and $mbcd(b^2-4c)$
nonzero.

\ \ {\tt Type I}. $\sum_{k=0}^\infty(a+dk)\bi{2k}k^2T_k(b,c)/m^k$.

\ \ {\tt Type II}.
$\sum_{k=0}^\infty(a+dk)\bi{2k}k\bi{3k}kT_k(b,c)/m^k$.

\ \ {\tt Type III}.
$\sum_{k=0}^\infty(a+dk)\bi{4k}{2k}\bi{2k}kT_k(b,c)/m^k$.

\ \ {\tt Type IV}.
$\sum_{k=0}^\infty(a+dk)\bi{2k}{k}^2T_{2k}(b,c)/m^k$.

\ \ {\tt Type V}.
$\sum_{k=0}^\infty(a+dk)\bi{2k}{k}\bi{3k}kT_{3k}(b,c)/m^k$.
\medskip

During October 1-3, 2011, I introduced two new kinds of series for $1/\pi$:
\medskip
\ \ {\tt Type VI}.
$\sum_{k=0}^\infty(a+dk)T_{k}(b,c)^3/m^k,$
\medskip
\ \ {\tt Type VII}.
$\sum_{k=0}^\infty(a+dk)\bi{2k}kT_{k}(b,c)^2/m^k,$
\medskip
\noindent where $a,b,c,d,m$ are integers and $mbcd(b^2-4c)$ is nonzero.

Recall that a series $\sum_{k=0}^\infty a_k$ is said to converge at a geometric rate with ratio $r$
if $\lim_{k\to+\infty}a_{k+1}/a_k=r\in(0,1)$.

\proclaim{Conjecture I {\rm (Z. W. Sun [S14b])}} We have the
following identities:
$$\align\sum_{k=0}^\infty\f{30k+7}{(-256)^k}\bi{2k}k^2T_k(1,16)=&\f{24}{\pi},\tag I1
\\\sum_{k=0}^\infty\f{30k+7}{(-1024)^k}\bi{2k}k^2T_k(34,1)=&\f{12}{\pi},\tag $\underline{\t{I2}}$
\\\sum_{k=0}^\infty\f{30k-1}{4096^k}\bi{2k}k^2T_k(194,1)=&\f{80}{\pi},\tag I3
\\\sum_{k=0}^\infty\f{42k+5}{4096^k}\bi{2k}k^2T_k(62,1)=&\f{16\sqrt3}{\pi}.\tag $\underline {\t{I4}}$
\endalign$$
\endproclaim
\Remark. The series (I1)-(I4) converge at geometric rates with
ratios $-9/16$, $-9/16$, $49/64$, $1/4$ respectively.

 \proclaim {Conjecture II {\rm (Z. W. Sun [S14b])}} We have
$$\align\sum_{k=0}^\infty\f{15k+2}{972^k}\bi{2k}k\bi{3k}k
T_k(18,6)=&\f{45\sqrt3}{4\pi},\tag $\underline{\t{II1}}$
\\\sum_{k=0}^\infty\f{91k+12}{10^{3k}}\bi{2k}k\bi{3k}kT_k(10,1)=&\f{75\sqrt3}{2\pi},\tag
II2
\\\sum_{k=0}^\infty\f{15k-4}{18^{3k}}\bi{2k}k\bi{3k}kT_k(198,1)=&\f{135\sqrt3}{2\pi},\tag
II3
\\\sum_{k=0}^\infty\f{42k-41}{30^{3k}}\bi{2k}k\bi{3k}kT_k(970,1)=&\f{525\sqrt3}{\pi},\tag
II4
\\\sum_{k=0}^\infty\f{18k+1}{30^{3k}}\bi{2k}k\bi{3k}kT_k(730,729)=&\f{25\sqrt3}{\pi},\tag II5
\endalign$$
$$\align\sum_{k=0}^\infty\f{6930k+559}{102^{3k}}\bi{2k}k\bi{3k}kT_k(102,1)=&\f{1445\sqrt6}{2\pi},\tag
II6
\\\sum_{k=0}^\infty\f{222105k+15724}{198^{3k}}\bi{2k}k\bi{3k}kT_k(198,1)=&\f{114345\sqrt3}{4\pi},\tag II7
\\\sum_{k=0}^\infty\f{390k-3967}{102^{3k}}\bi{2k}k\bi{3k}kT_k(39202,1)=&\f{56355\sqrt3}{\pi},\tag II8
\\\sum_{k=0}^\infty\f{210k-7157}{198^{3k}}\bi{2k}k\bi{3k}kT_k(287298,1)=&\f{114345\sqrt{3}}{\pi},\tag II9
\endalign$$
and
$$\align\sum_{k=0}^\infty\f{45k+7}{24^{3k}}\bi{2k}k\bi{3k}kT_k(26,729)=&\f8{3\pi}(3\sqrt3+\sqrt{15}),\tag II10
\\\sum_{k=0}^\infty\f{9k+2}{(-5400)^k}\bi{2k}k\bi{3k}kT_k(70,3645)=&\f{15\sqrt3+\sqrt{15}}{6\pi},\tag $\underline{\t{II11}}$
\\\sum_{k=0}^\infty\f{63k+11}{(-13500)^k}\bi{2k}k\bi{3k}kT_k(40,1458)=&\f{25}{12\pi}(3\sqrt3+4\sqrt6),\tag II12
\endalign$$
\endproclaim
\Remark. The series (II1)-(II12) converge at geometric rates with
ratios
$$\gather\f{9+\sqrt6}{18},\ \f{81}{250},\ \f{25}{27},\ \f{243}{250},\ \f{98}{125},
\ \f{13}{4913},\ \f{25}{35937},
\\\f{9801}{9826},\ \f{71825}{71874},\ \f5{32},\ -\f{35+27\sqrt5}{100},
\ -\f{20+27\sqrt2}{250}\endgather$$
respectively.

\proclaim{Conjecture III {\rm (Z. W. Sun [S14b])}} We have the following formulae:
$$\align\sum_{k=0}^\infty\f{85k+2}{66^{2k}}\bi{4k}{2k}\bi{2k}kT_k(52,1)=&\f{33\sqrt{33}}{\pi},\tag III1
\\\sum_{k=0}^\infty\f{28k+5}{(-96^2)^k}\bi{4k}{2k}\bi{2k}kT_k(110,1)=&\f{3\sqrt6}{\pi},\tag III2
\\\sum_{k=0}^\infty\f{40k+3}{112^{2k}}\bi{4k}{2k}\bi{2k}kT_k(98,1)=&\f{70\sqrt{21}}{9\pi},\tag $\underline{\t{III3}}$
\\\sum_{k=0}^\infty\f{80k+9}{264^{2k}}\bi{4k}{2k}\bi{2k}kT_k(257,256)=&\f{11\sqrt{66}}{2\pi},\tag III4
\\\sum_{k=0}^\infty\f{80k+13}{(-168^2)^k}\bi{4k}{2k}\bi{2k}kT_k(7,4096)=&\f{14\sqrt{210}+21\sqrt{42}}{8\pi},\tag
$\underline{\t{III5}}$
\endalign$$
and
$$\align\sum_{k=0}^\infty\f{760k+71}{336^{2k}}\bi{4k}{2k}\bi{2k}kT_k(322,1)=&\f{126\sqrt{7}}{\pi},\tag III6
\\\sum_{k=0}^\infty\f{10k-1}{336^{2k}}\bi{4k}{2k}\bi{2k}kT_k(1442,1)=&\f{7\sqrt{210}}{4\pi},\tag III7
\\\sum_{k=0}^\infty\f{770k+69}{912^{2k}}\bi{4k}{2k}\bi{2k}kT_k(898,1)=&\f{95\sqrt{114}}{4\pi},\tag III8
\\\sum_{k=0}^\infty\f{280k-139}{912^{2k}}\bi{4k}{2k}\bi{2k}kT_k(12098,1)=&\f{95\sqrt{399}}{\pi},\tag III9
\\\sum_{k=0}^\infty\f{84370k+6011}{10416^{2k}}\bi{4k}{2k}\bi{2k}kT_k(10402,1)=&\f{3689\sqrt{434}}{4\pi},\tag III10
\\\sum_{k=0}^\infty\f{8840k-50087}{10416^{2k}}\bi{4k}{2k}\bi{2k}kT_k(1684802,1)=&\f{7378\sqrt{8463}}{\pi},\tag III11
\\\sum_{k=0}^\infty\f{11657240k+732103}{39216^{2k}}\bi{4k}{2k}\bi{2k}kT_k(39202,1)=&\f{80883\sqrt{817}}{\pi},\tag III12
\\\sum_{k=0}^\infty\f{3080k-58871}{39216^{2k}}\bi{4k}{2k}\bi{2k}kT_k(23990402,1)=&\f{17974\sqrt{2451}}{\pi}.\tag III13
\endalign$$
\endproclaim
\Remark. The series (III1)-(III13) converge at geometric rates with
ratios
$$\gather\f{96}{121},\ -\f 79,\ \f{25}{49},\ \f{289}{1089},\ -\f{15}{49},\ \f 9{49},\ \f{361}{441},
\\\f{25}{361},\ \f{3025}{3249},\ \f{289}{47089},\ \f{421201}{423801},\  \f{1089}{667489},\ \f{5997601}{6007401}
\endgather$$
respectively. I thank Prof. Qing-Hu Hou (at Nankai University) for helping me check (III9) numerically.

\proclaim{Conjecture IV {\rm (Z. W. Sun [S14b])}} We have
$$\align
\sum_{k=0}^\infty\f{26k+5}{(-48^2)^k}\bi{2k}k^2T_{2k}(7,1)=&\f{48}{5\pi},\tag IV1
\\\sum_{k=0}^\infty\f{340k+59}{(-480^2)^k}\bi{2k}k^2T_{2k}(62,1)=&\f{120}{\pi},\tag IV2
\\\sum_{k=0}^\infty\f{13940k+1559}{(-5760^2)^k}\bi{2k}k^2T_{2k}(322,1)=&\f{4320}{\pi},\tag IV3
\\\sum_{k=0}^\infty\f{8k+1}{96^{2k}}\bi{2k}k^2T_{2k}(10,1)=&\f{10\sqrt{2}}{3\pi},\tag $\underline{\t{IV4}}$
\\\sum_{k=0}^\infty\f{10k+1}{240^{2k}}\bi{2k}k^2T_{2k}(38,1)=&\f{15\sqrt6}{4\pi},\tag IV5
\\\sum_{k=0}^\infty\f{14280k+899}{39200^{2k}}\bi{2k}k^2T_{2k}(198,1)=&\f{1155\sqrt{6}}{\pi},\tag IV6
\\\sum_{k=0}^\infty\f{120k+13}{320^{2k}}\bi{2k}k^2T_{2k}(18,1)=&\f{12\sqrt{15}}{\pi},\tag IV7
\\\sum_{k=0}^\infty\f{21k+2}{896^{2k}}\bi{2k}k^2T_{2k}(30,1)=&\f{5\sqrt7}{2\pi},\tag IV8
\\\sum_{k=0}^\infty\f{56k+3}{24^{4k}}\bi{2k}k^2T_{2k}(110,1)=&\f{30\sqrt7}{\pi},\tag IV9
\\\sum_{k=0}^\infty\f{56k+5}{48^{4k}}\bi{2k}k^2T_{2k}(322,1)=&\f{72\sqrt7}{5\pi},\tag IV10
\\\sum_{k=0}^\infty\f{10k+1}{2800^{2k}}\bi{2k}k^2T_{2k}(198,1)=&\f{25\sqrt{14}}{24\pi},\tag IV11
\\\sum_{k=0}^\infty\f{195k+14}{10400^{2k}}\bi{2k}k^2T_{2k}(102,1)=&\f{85\sqrt{39}}{12\pi},\tag  IV12
\endalign$$
$$\align\sum_{k=0}^\infty\f{3230k+263}{46800^{2k}}\bi{2k}k^2T_{2k}(1298,1)=&\f{675\sqrt{26}}{4\pi},\tag IV13
\\\sum_{k=0}^\infty\f{520k-111}{5616^{2k}}\bi{2k}k^2T_{2k}(1298,1)=&\f{1326\sqrt3}{\pi},\tag IV14
\\\sum_{k=0}^\infty\f{280k-149}{20400^{2k}}\bi{2k}k^2T_{2k}(4898,1)=&\f{330\sqrt{51}}{\pi},\tag IV15
\\\sum_{k=0}^\infty\f{78k-1}{28880^{2k}}\bi{2k}k^2T_{2k}(5778,1)=&\f{741\sqrt{10}}{20\pi},\tag IV16
\\\sum_{k=0}^\infty\f{57720k+3967}{439280^{2k}}\bi{2k}k^2T_{2k}(5778,1)=&\f{2890\sqrt{19}}{\pi},\tag IV17
\\\sum_{k=0}^\infty\f{1615k-314}{243360^{2k}}\bi{2k}k^2T_{2k}(54758,1)=&\f{1989\sqrt{95}}{4\pi},\tag IV18
\\\sum_{k=0}^\infty\f{34k+5}{4608^k}\bi{2k}k^2T_{2k}(10,-2)=&\f{12\sqrt6}{\pi},\tag IV19
\\\sum_{k=0}^\infty\f{130k+1}{1161216^k}\bi{2k}k^2T_{2k}(238,-14)=&\f{288\sqrt2}{\pi},\tag IV20
\\\sum_{k=0}^\infty\f{2380k+299}{(-16629048064)^k}\bi{2k}k^2T_{2k}(9918,-19)=&\f{860\sqrt7}{3\pi}.\tag IV21
\endalign$$
\endproclaim
\Remark. The series (IV1)-(IV21) converge at geometric rates with
ratios
$$\gather-\f 9{16},\ -\f{64}{225},\ -\f{81}{1600},\ \f{1}{4},\ \f4 9,\ \f{1}{2401},\ \f1{16},\ \f1{49},\ \f{49}{81},\ \f{81}{256},\ \f{4}{49},\ \f1{625},
\\ \f1{81},\ \f{625}{729},\ \f{2401}{2601},\ \f{83521}{130321},\ \f1{361},
\ \f{1874161}{2313441},\ \f{3}8,\ \f{25}{32},\ -\f{175}{1849}\endgather$$
respectively. I conjecture that (IV1)-(IV18) have exhausted all identities of the form
$$\sum_{k=0}^\infty(a+dk)\f{\bi{2k}k^2T_{2k}(b,1)}{m^k}=\f C{\pi}$$
with $a,d,m\in\Z$, $b\in\{1,3,4,\ldots\}$, $d>0$, and $C^2$ positive and rational.

\proclaim{Conjecture V {\rm (Z. W. Sun [S14b])}} We have the formula
$$\sum_{k=0}^\infty\f{1638k+277}{(-240)^{3k}}\bi{2k}k\bi{3k}kT_{3k}(62,1)=\f{44\sqrt{105}}{\pi}.\tag V1$$
\endproclaim
\Remark. The series (V1) converges at a geometric rate with ratio $-64/125$.
\medskip

Note that [CWZ1] contains complete proofs of (I2), (I4), (II1), (II11), (III3) and (III5).
Also, a detailed proof of (IV4) was given in [WZ]. The most crucial parts of such proofs involve modular equations,
so (in my opinion) a complete proof should contain all the details involving modular equation.

\proclaim{Conjecture VI {\rm (Z. W. Sun [S14b])}} We have the following formulae:
$$\align
\sum_{k=0}^\infty\f{66k+17}{(2^{11}3^3)^{k}}T_k^3(10,11^2)=&\f{540\sqrt2}{11\pi},\tag VI1
\\\sum_{k=0}^\infty\f{126k+31}{(-80)^{3k}}T_k^3(22,21^2)=&\f{880\sqrt5}{21\pi},\tag VI2
\\\sum_{k=0}^\infty\f{3990k+1147}{(-288)^{3k}}T_k^3(62,95^2)=&\f{432}{95\pi}(195\sqrt{14}+94\sqrt2).\tag VI3
\endalign$$
\endproclaim
\Remark. The series (VI1)-(VI3) converge at geometric rates with ratios
$$\f{16}{27},\ -\f{64}{125},\ -\f{343}{512}$$
respectively. I would like to offer \$300 as the prize for the person (not joint authors)
who can provide first rigorous proofs of all the three identities (VI1)-(VI3).

\proclaim{Conjecture VII {\rm (Z. W. Sun [S14b])}} We have the following formulae:
$$\align
\sum_{k=0}^\infty\f{221k+28}{450^{k}}\bi{2k}kT_k^2(6,2)=&\f{2700}{7\pi},\tag$\underline{\t{VII1}}$
\\\sum_{k=0}^\infty\f{24k+5}{28^{2k}}\bi{2k}kT_k^2(4,9)=&\f{49}{9\pi}(\sqrt3+\sqrt6),\tag VII2
\\\sum_{k=0}^\infty\f{560k+71}{22^{2k}}\bi{2k}kT_k^2(5,1)=&\f{605\sqrt7}{3\pi},\tag$\underline{\t{VII3}}$
\\\sum_{k=0}^\infty\f{3696k+445}{46^{2k}}\bi{2k}kT_k^2(7,1)=&\f{1587\sqrt7}{2\pi},\tag$\underline{\t{VII4}}$
\\\sum_{k=0}^\infty\f{56k+19}{(-108)^k}\bi{2k}kT_k^2(3,-3)=&\f{9\sqrt7}{\pi},\tag$\underline{\t{VII5}}$
\endalign$$
$$\align\sum_{k=0}^\infty\f{450296k+53323}{(-5177196)^k}\bi{2k}kT_k^2(171,-171)=&\f{113535\sqrt7}{2\pi},\tag$\underline{\t{VII6}}$
\\\sum_{k=0}^\infty\f{2800512k+435257}{434^{2k}}\bi{2k}kT_k^2(73,576)=&\f{10406669}{2\sqrt6\,\pi}.\tag VII7
\endalign$$
\endproclaim
\Remark. The series (VII1)-(VII7) converge at geometric rates with ratios
$$\f{88+48\sqrt2}{225},\ \f{25}{49},\ \f{49}{121},\ \f{81}{529},\ -\f 79,\ -\f{175}{7569},\ \f{14641}{47089}.$$
respectively. W. Zudilin [Zu] discussed (VII1) and (VII2)-(VII6) with the help of S. Cooper's work [Co].

\heading{4. Historical notes on the 61 series in Section
3}\endheading

 I discovered most of those conjectural series for $1/\pi$ in Section 3 during Jan. and Feb. in 2011.
 Series of type VI and VII were introduced in October 2011.
 All my conjectural series in Section 2 came from a combination of my philosophy, intuition, inspiration, experience and computation.

 In the evening of Jan. 1, 2011 I figured out the asymptotic
 behavior of $T_n(b,c)$ with $b$ and $c$ positive. (Few days later I learned the Laplace-Heine asymptotic formula for
 Legendre polynomials and hence knew that my conjectural main term of $T_n(b,c)$ as $n\to+\infty$ is indeed correct.)

 The story of new series for $1/\pi$ began with (I1) which was found in the early morning of Jan. 2, 2011 immediately after I waked up on the bed. On
 Jan 4 I announced this via a message to {\tt Number Theory Mailing List} as well as the initial version of [S14b] posted to {\tt arXiv}.
 In the subsequent two weeks I communicated with some experts on
 $\pi$-series and wanted to know whether they could prove my
 conjectural (I1).  On Jan. 20, it seemed clear that series like
 (I1) could not be easily proved by the current known methods used to establish
 Ramanujan-type series for $1/\pi$.

 Then, I discovered (II1) on Jan. 21 and (III3) on Jan. 29. On Feb. 2 I
 found (IV1) and (IV4). Then, I discovered (IV2) on Feb. 5. When I waked up in the early morning of Feb. 6,
 I suddenly realized a (conjectural) criterion for the existence of series for $1/\pi$ of type IV with $c=1$.
 Based on this criterion,  I found (IV3), (IV5)-(IV10) and (IV12) on
 Feb. 6, (IV11) on Feb. 7, (IV13) on Feb. 8,  (IV14)-(IV16) on Feb. 9, and (IV17) on Feb. 10.
 On Feb. 14 I discovered (I2)-(I4) and (III4). I found the
 sophisticated (III5) on Feb. 15. As for series of type IV,
 I discovered the largest example (IV18) on Feb. 16., and conjectured
 that the 18 series in Conj. IV have exhausted all those series for $1/\pi$ of type IV with $c=1$.
 On Feb. 18 I found (II2), (II5)-(II7), (II10) and (II12).

 On Feb. 21 I informed many experts on $\pi$-series (including Gert Almkvist) my list of the 34 conjectural series for
$1/\pi$ of types I-IV and predicted that there are totally about 40 such series. On Feb. 22 I found (II11) and
 (II3)-(II4); on the same day, motivated by my conjectural (II2),(II5)-(II7), (II10) and (II12) discovered on Feb 18,
 G. Almkvist found the following two series of type II that I missed:
 $$\sum_{k=0}^\infty\f{42k+5}{18^{3k}}\bi{2k}k\bi{3k}kT_k(18,1)=\f{54\sqrt3}{5\pi}\tag A1$$
 and
 $$\sum_{k=0}^\infty\f{66k+7}{30^{3k}}\bi{2k}k\bi{3k}kT_k(30,1)=\f{50\sqrt2}{3\pi}.\tag A2$$
On Feb. 22, Almkvist also pointed out that my conjectural identity
(II2) can be used to compute an arbitrary decimal digit of
$\sqrt3/\pi$ without computing the earlier digits.

 On Feb. 23 I discovered (V1), which is the unique example of series for $1/\pi$ of type V that I can find.

 On Feb. 25 and Feb. 26, I found (II8) and (II9) respectively. These two series converge very slowly.

 On August 11, I discovered (III6)-(III9) that I missed during Jan.-Feb.
 (III10)-(III13) were found by me on Sept. 21, 2011. Note that (III13) converges very slow.

 On Oct. 1,  I discovered (VI2) and (VI3), then I found (VI1) on the next day.

 I figured out (VII1)-(VII4), (VII5) and (VII6) on Oct. 3, 4 and 5 respectively.
 On Oct. 13, 2011 I discovered (VII7).

 On Oct. 16 James Wan informed me the preprints [CWZ1] and [WZ] on my conjectural series of types I-V.
 I admit that these two papers contain complete proofs of (I2), (I4), (II1), (II11), (III3), (III5) and (IV4).
 Note also that [CWZ2] was motivated by the authors' study of my conjectural (III5).

 On Oct. 7, 2012 I found (IV19)-(IV21) which involve $T_{2k}(b,c)$ with $c<0$.

 My paper [S14b] containing the 61 series in Section 3 was finally published in 2014.

\heading{5. A technique for producing more series for $1/\pi$}\endheading

For a sequence $a_0,a_1,a_2,\ldots$ of complex numbers, define
$$a_n^*=\sum_{k=0}^n\bi nk(-1)^ka_k\quad\t{for all}\ n\in\N=\{0,1,2,\ldots\}$$
and call $\{a_n^*\}_{n\in\N}$ the {\it dual sequence} of $\{a_n\}_{n\in\N}.$
It is well known that $(a_n^*)^*=a_n$ for all $n\in\N$.

There are many series for $1/\pi$ of the form
$$\sum_{k=0}^\infty (bk+c)\f{\bi{2k}ka_k}{m^k}=\f {C}{\pi},$$
where $a_k,b,c,C$ and $m\not=0$ are real numbers
(see Section 2-3 for many such series). On March 10, 2011, I realized that if $|m-4|>4$ then
$$\sum_{n=0}^\infty(bmn+2b+(m-4)c)\f{\bi{2n}na_n^*}{(4-m)^n}
=(m-4)\sqrt{\f{m-4}m}\sum_{k=0}^\infty(bk+c)\f{\bi{2k}ka_k}{m^k}.\tag5.1$$
(For the reason, see [S14c, Section 1].)
Thus, if $m>8$ or $m<0$ then
$$\aligned &\sum_{k=0}^\infty(bk+c)\f{\bi{2k}ka_k}{m^k}=\f C{\pi}
\\\Longrightarrow &\sum_{k=0}^\infty(bmk+2b+(m-4)c)\f{\bi{2k}ka_k^*}{(4-m)^k}=\f{(m-4)C}{\pi}\sqrt{\f{m-4}m}.
\endaligned\tag5.2$$
\medskip

 {\it Example} 5.1. Let $a_n=\bi{2n}nT_n(1,16)$ for all $n\in\N$. Then
$$a_n^*=\sum_{k=0}^n\bi nk\bi{2k}k(-1)^kT_k(1,16)\quad\t{for}\ n=0,1,2,\ldots.$$
Thus, by (5.2), the identity (I1) in Section 3 implies that
$$\sum_{k=0}^\infty(48k+11)\f{\bi{2k}ka_k^*}{260^k}=\f{39\sqrt{65}}{8\pi}.$$

\enddocument
\end